\documentclass[11pt,twoside,letter]{article}

\usepackage{color}
\usepackage{float}
\usepackage{amsmath,amssymb}
\usepackage{mathrsfs}
\usepackage{enumerate}
\usepackage{srcltx}
\usepackage{graphicx}

\usepackage{color}

\newcommand {\Gc}  {\mathcal{G}}
\newcommand {\Lc}  {\mathcal{L}}

\definecolor{red}{rgb}{0.9,0,0}
\definecolor{blue}{rgb}{0.2,0.2,0.7}
\definecolor{green}{rgb}{0.0,0.5,0.2}
\definecolor{darkblue}{rgb}{0.2,0.2,0.5}
\definecolor{orange}{rgb}{1,0.5,0}

\begin{document}

\title{\vspace{-3cm}
\LARGE \bf
A Class of\\ Minimal Generically Universally Rigid Frameworks
}

\author{Scott D. Kelly$^{\natural}$ and Andrea Micheletti$^{\sharp}$
}

\date{\small\today}

\maketitle

\vspace{-1.0cm}
\begin{center}
{\small

$^\natural$ Department of Mechanical Engineering and Engineering Sciences\\
University of North Carolina at Charlotte, USA\\
{\tt\small scott at kellyfish.net}\\[10pt]

$^\sharp$ Dipartimento di Ingegneria Civile e Ingegneria Informatica\\
Universit\`a di Roma Tor Vergata, Italy\\
{\tt\small micheletti at ing.uniroma2.it}

}
\end{center}

\pagestyle{myheadings}
\markboth{S.~D.~Kelly and A.~Micheletti}
{A class of minimal generically universally rigid frameworks}

\section*{Abstract}
Following a review of related results in rigidity theory, we provide a construction to obtain generically universally rigid frameworks with the minimum number of edges, for any given set of $n$ nodes in two or three dimensions. When a set of edge-lengths is compatible with only one configuration in $d$-dimensions, the framework is globally rigid. When that configuration is unique even if embedded in a higher dimensional space, the framework is universally rigid. In case of generic configurations, where the nodal coordinates are algebraically independent
, the minimum number of edges required is equal to $dn-d(d+1)/2+1$, that is, $2n-2$ for $d=2$, and $3n-5$ for $d=3$. Our contribution is a specific construction for this case by introducing a class of frameworks generalizing that of Gr\"{u}nbaum polygons. The construction applies also to nongeneric configurations, although in this case the number of edges is not necessarily the minimum. One straightforward application is the design of wireless sensor networks or multi-agent systems with the minimum number of communication links.\\[5pt]
\noindent {\em Keywords: global rigidity, universal rigidity, Gr\"unbaum polygon, sensor network, multi-agent system}


\section{Introductory definitions and results}

The analysis of molecular conformations, the localization of a wireless sensor network, or the coordination of a multi-agent system are instances of problems where the rigidity of the underlying framework plays an important role (see, e.g., \cite{Jacobs1999,Jackson2009,Zhu2010,Smith2009,Nabet2009}). The simplest and most often adopted model is that of bar-and-joint framework, i.e. a set of points in space together with a set of distance assignments. In this section we review some relevant definitions and results. The following section is dedicated to our main result.

Rigorously, a framework is a graph together with a configuration. Let $\Gc$ be a graph on $n$ nodes with $e$ edges. The edge connecting the $i$-th node with the $j$-th node is denoted by $ij$, considering only finite and undirected graphs, without loops or multiple edges. A configuration assigns to every node a point in the $d$-dimensional Euclidean space $E^d$.
Let $p$ be the vector collecting all nodal coordinates with respect to a given reference frame. Then $(\Gc,p)$ is the framework with graph $\Gc$ and configuration $p$. Associated to a framework is the set $\Lc(\Gc,p)$ of the squared edge lengths.

\begin{figure}[h!]
\begin{center}
\includegraphics[scale=1.83]{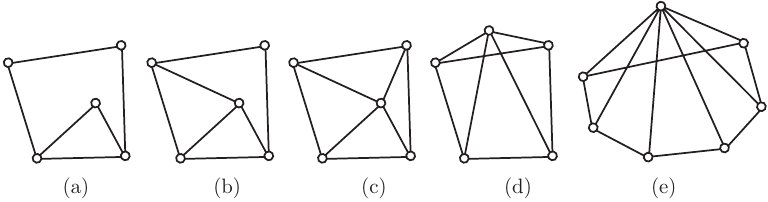}    
\caption{Frameworks in two dimensions belonging to different rigidity classes: flexible (a), rigid (b), globally rigid (c), and universally rigid (d,e). The last two frameworks are known as Gr\"unbaum polygons \cite{Grunbaum1978}.}
\label{f1}
\end{center}
\end{figure}

Different classes of rigididy can be defined for frameworks (Fig.~\ref{f1}). Here, we mostly follow the treatment in \cite{Connelly2005,Connelly2009c}. A configuration $q$ is {\em admissible} for $(\Gc,p)$ if $\Lc(G,q)=\Lc(G,p)$.
Two configurations $p$ and $q$ are {\em congruent}, and we write $p\equiv q$, if in both configurations all distances between node pairs are the same.
Equivalently, two configurations are congruent if they differ by an isometry of $E^d$, i.e a composition of translations, rotations and reflections.
A framework $(\Gc,p)$ is {\em rigid} if
any admissible configuration $q$ which is close to $p$
is congruent to $p$.

The Jacobian of $\Lc(\Gc,p)$ with respect to $p$ is the ($e$-by-$dn$) {\em rigidity matrix}, $R$. A framework is {\em infinitesimally rigid} if the rank of $R$ is equal to $nd-d(d+1)/2$, or equivalently, if the nullspace of $R$ contains only {\em rigid velocities}, i.e. nodal velocities in a rigid motion.
%
%
The $m$ independent vectors in the nullspace of $R$ which are not rigid velocities are called {\em flexes} (Fig.~\ref{f1} a), so that a framework is infinitesimally rigid if it has no flexes ($m=0$, Fig.~\ref{f1} b).
A framework $(\Gc,p)$ is {\em globally rigid} if any admissible configuration $q$ is congruent to $p$ (Fig.~\ref{f1} c). A framework is {\em universally rigid} if it is globally rigid in all dimensions (Fig.~\ref{f1} d, e). 


For each rigidity class, the generic and nongeneric cases can be distinguished.
A configuration is {\em generic} if the coordinates in $p$ are algebraically independent, i.e. if the nodal coordinates do not satisfy any nontrivial polynomial equation with integer coefficients.
In nongeneric configurations, rigidity properties are more difficult to predict. Standard examples are those where a flexible framework becomes rigid, or where a framework loses global rigidity, when passing from a generic configuration to a nongeneric one \cite{Connelly2009c}.

A framework $(\Gc,p)$ is {\em generically rigid} if it is rigid and $p$ is generic. Generic rigidity is a property of the graph, not the configuration.
The minimum number of edges necessary for generic rigidity are $2n-3$ in 2D and $3n-6$ in 3D. Intuitively, in 2D, we can start with an edge connecting two nodes, then iteratively adding one node connected to the other nodes by two noncollinear edges. In 3D, we can start with a nondegenerate {\em triangle} (three vertices and three edges), then iteratively adding one {\em tripod}, i.e. a node connected to the other nodes by three noncoplanar edges. These constructions constitute particular {\em Henneberg sequences} \cite{Henneberg1911,Eren2004}: sequences of operations which preserve minimal generic rigidity.


The characterization of global and universal rigidity has been given in the literature in terms of {\em stress}. A stress $\omega$ is an assignment of a real number
to each edge
of the framework. This number is interpreted as the axial force carried by the edge, divided by the length of that edge. A {\em selfstress} for $(\Gc,p)$ is a stress satisfying the equilibrium conditions at each node in absence of external forces: the composition of the forces concurrent at a node is null, for every node.
It is easy to see that selfstresses belong to the nullspace of $R^T$, also called the equilibrium matrix. It follows that the number of independent selfstresses $s$ and mechanisms $m$ are related to $n$ and $e$ by the following rule
\begin{equation}\label{Maxwellrule}
dn-d(d+1)/2-e=m-s\,,
\end{equation}
where $d(d+1)/2$ is the number of independent rigid motions in $E^d$. This rule follows from the orthogonality of the fundamental subspaces of $R$
\cite{pca86}.

A fundamental object is the ($n$-by-$n$) {\em stress matrix}, $\Omega$, that is, the weighted Laplacian of the graph, with weights given by the selfstress values on the edges. Since weights can be either positive or negative, many standard results on positively-weighted Laplacians do not apply.
A first characterization of universal rigidity has been given by Connelly through $\Omega$ by introducing the notion of superstability \cite{Connelly1982, Connelly2009c}.
A framework in $E^d$ with the affine span of the nodal position vectors
being all $E^d$ and a nonzero selfstress is {\em superstable} if: 1) $\Omega$ is positive semidefinite; 2) $\Omega$ has rank $n-d-1$; 3) there are no affine admissible motions. Connelly has shown that {\em every superstable framework is universally rigid}.
More recently, the converse statement has been shown in the generic case \cite{Gortler2010a}: {\em every  universally rigid framework $(\Gc,p)$ with $p$ generic and $n\geq d+2$ is superstable}.
The minimum number of edges necessary for universal rigidity is then easy to find. First consider a {\em simplex} in $E^d$, that is, a framework on the complete graph on $d+1$ nodes, e.g. triangles in $E^2$ or tetrahedra in $E^3$. Simplices, and all frameworks on complete graphs, are universally rigid by definition, since admissible configurations must be congruent to each other. Every generic universally rigid framework in $E^d$ which is not a simplex, i.e. it has at least $d+2$ nodes, admits at least one independent selfstress, $s\geq 0$, to have a non-null $\Omega$, and also it has no flexes, $m=0$, since $p$ is generic. It follows from (\ref{Maxwellrule}) that in a generically universally rigid framework the number of edges satisfies
\begin{equation}\label{edgesnumber}
e\geq dn-d(d+1)/2+1\,,
\end{equation}
i.e. $e\geq 2n-2$ in $d=2$ or $e\geq 3n-5$ in $d=3$.

A less strict condition for universal rigidity consists in requiring a configuration to be {\em general}. A configuration in $E^d$ is general if no $d+1$ nodes are affinely dependent, e.g. there are no three collinear nodes if $d=2$, or there are no three collinear nodes and no four coplanar nodes if $d=3$. In this case we have that \cite{Alfakih2010}: {\em a framework $(\Gc,p)$ with $p$ general and $n\geq d+2$ is universally rigid if there is a nonzero selfstress whose stress matrix is positive semi-definite with rank $n-d-1$}.
It has been shown in \cite{Alfakih2011} that the converse of this theorem holds for $(d+1)$-lateration graphs, i.e. graphs obtained from a simplex by applying a sequence of $(d+1)$-valent node additions, i.e. the addition of a node connected by $d+1$ edges to the other nodes. Notice that the number of edges of frameworks obtained in this way is $e=(d+1)n-(d+2)(d+1)/2$, that is $e=3n-6$ for $d=2$ and $e=4n-10$ for $d=3$. Notice that, for  large $n$, these values of $e$ are $50\%$ and $33\%$ higher than the minimum value given by \eqref{edgesnumber}, respectively for $d=2$ and $d=3$.

Known universally rigid frameworks are: cablenets and frameworks in the shape of convex polygons \cite{Connelly1982}; frameworks in the shape of centrally symmetric polyhedra \cite{Lovasz2001,Bezdek2006}. It has also been shown that given two universally rigid frameworks, it is possible to combine them into a universally rigid assembly if they have $d+1$ nodes in common \cite{Ratmanski2010}.

In the next section we show, without relying on existing results,
that it is always possible to construct frameworks on $n$ given nodes in $E^2$ or $E^3$ with the minimum number of edges \eqref{edgesnumber}, irrespective of the generic/nongeneric property of the configuration. Such frameworks belong to a new class which generalize that of Gr\"unbaum polygons \cite{Grunbaum1978}. One straighforward application of this construction is the realization of sensor networks or multi-agent systems with the minimum number of communication links.

\begin{figure}[h!]
\begin{center}
\includegraphics[scale=1.83]{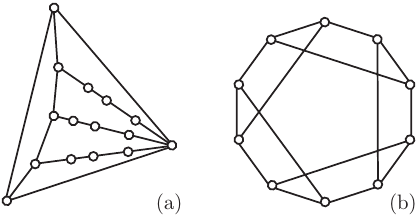}    
\caption{Two examples of nongeneric minimal systems. These frameworks have less than $2n-2$ edges and are universally rigid in two dimensions. In (a) a cable net is attached to a simplex; in (b) the outer edges form a regular polygon.}
\label{f2b}
\end{center}
\end{figure}

In case of nongeneric configurations, it is possible to decrease the number of edges further (Fig.~\ref{f2b}); however, limited theoretical results are available for nongeneric systems. In a follow-up paper we plan to give a list of possible constructions for this case.
%



\section{Construction of minimal generic universally rigid frameworks}\label{results}

\subsection{Nonconvex Gr\"unbaum polygons}

\begin{figure}[h!]
\begin{center}
\includegraphics[scale=1.83]{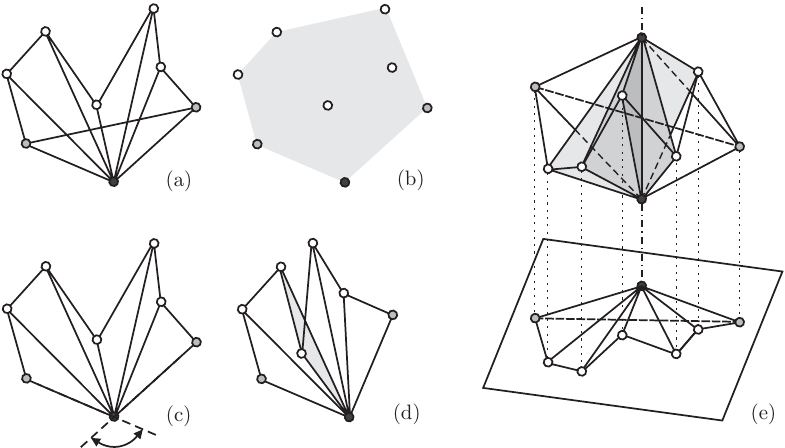}    
\caption{A nonconvex Gr\"unbaum polygon is shown in (a). Frameworks of this kind can be obtained by starting from a set of nodes in 2D, then by finding three consecutive nodes on the boundary of the convex hull of the set (b), the middle node is the {\em center} (in black), the other two nodes are its {\em neighbors} (in grey). In (c) all the nodes are connected to the center, and edges added to form adjacent triangles, so as to obtain a {\em fan}. In (d) one of the admissible configurations of this fan is shown, obtained by flipping the triangle shown (in light grey). In the completely unfolded fan (c) the angle shown is maximized. The universally rigid framework in (a) is obtained by adding the last edge between the two neighbors of the center. The three dimensional framework in (e) can be obtained in analogous way (see description in the text); the end nodes of the {\em central edge} are shown in black, and the neighbors in grey. A nonconvex Gr\"unbaum polygon can be obtained as a projection onto a plane along the direction of the central edge.}
\label{grun2d}
\end{center}
\end{figure}

Gr\"{u}nbaum polygons are frameworks obtained by placing nodes and edges respectively at the vertices and the sides of a convex polygon, then by choosing one node, the {\em center node} (in black in Fig.\,\ref{grun2d}\,a), and by connecting all the other nodes to it with an edge. The construction is completed by adding one edge connecting the two nodes neighboring the center.

We provide here a similar construction to assign $(2n-2)$ edges to any given a set of nodes in $E^2$ in order to obtain {\em nonconvex Gr\"unbaum polygons} (Fig.\,\ref{grun2d}).
First, the convex hull of the nodes is constructed and three consecutive vertices on its boundary coinciding with three nodes are chosen (Fig.\,\ref{grun2d}\,b). The middle one becomes the center to which all the other nodes are connected, with additional edges forming a contiguous sequence of triangles sharing the center as a vertex (Fig.\,\ref{grun2d}\,c). The last edge connects the two nodes neighboring the center (Fig.\,\ref{grun2d}\,a). We have the following result.

\

\noindent {\em Every nonconvex Gr\"unbaum polygon is universally rigid}.

\

\noindent {\em Proof.} Up to the addition of the last edge, the framework can be viewed as forming a {\em fan} which ``unfold'' from the center node (Fig.\,\ref{grun2d}\,c). This incomplete framework admits a number of configurations equal to $2^f$, where $f$ is the number of internal edges or {\em folds} of the fan (Fig.\,\ref{grun2d}\,d). The distance between the two neighboring nodes reaches a global maximum when the fan is completely unfolded. It follows that by adding the last edge between the two neighboring nodes, the unfolded configuration is unique.

By embedding this framework in a higher dimensional Euclidean space, the situation does not change. Since each triangle of a fan is universally rigid by itself and it can only rotate about a fold, relative to its neighboring triangles, the triangle inequality ensure that the distance between the two neighboring nodes has a global maximum only when the fan is flat, therefore the nonconvex Gr\"{u}mbaum polygon is universally rigid.  $_\Box$

Notice that this proof is valid for both convex and nonconvex Gr\"unbaum polygons. Notice also that the construction works even if the center is aligned with its neighbors, or if two or more fold are collinear. The result holds even if the configuration is nongeneric, the main requirement being that the center and its neighbors are on the boundary of the convex hull.

\subsection{Three-dimensional Gr\"unbaum frameworks}

In three dimensions we can obtain a perfectly analogous result for assigning $(3n-5)$ edges to a given a set of nodes in $E^3$. We construct the convex hull of this set. There will be at least four vertices of the hull forming two adjacent triangles, sharing one edge of the convex hull. The shared edge is the {\em central edge} of the framework, the two nodes on this edge are the {\em central nodes}, while the other two are the {\em neighboring nodes}. Now, we can add edges connecting each of the neighboring nodes to the central nodes. We do the same with the remaining nodes, by connecting them to the central nodes. In this way, we obtain a set of triangles in space, all sharing one edge. Then, for each couple of neighboring triangles, we add an edge between the nodes so as to form a tetrahedron. Finally, the last edge of this construction is added between the two neighboring nodes (Fig.\,\ref{grun2d} top right).

An easy way of visualizing this framework is to project it along the direction of the central edge onto a plane, resulting in a nonconvex Gr\"unbaum polygon (Fig.\,\ref{grun2d} right). In a way similar to what we have done before, we can consider the incomplete framework obtained by removing the last edge and argue that this admits a number of configuration equal to $2^f$, with $f$ defined for the projected framework as in the two-dimensional case. Among all this configurations, the one which is completely ``unfolded'' gives the maximum distance between the neighboring nodes, still using this term in analogy with the two-dimensional case. Once we add the last edge in this configuration, we obtain a globally rigid structure, which, by the triangle inequality is also universally rigid. We call frameworks obtained in this way {\em 3D Gr\"unbaum framework} and state the following result.

\

\noindent {\em Every 3D Gr\"unbaum framework is universally rigid}.

\

\begin{figure}[h!]
\begin{center}
\includegraphics[width=0.35\textwidth]{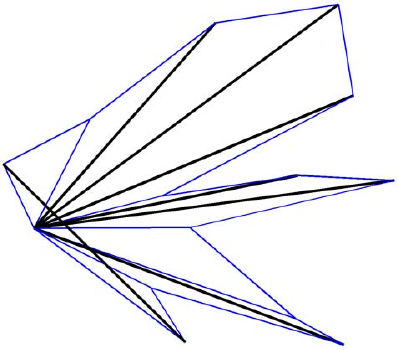}\,    
\quad\includegraphics[width=0.35\textwidth]{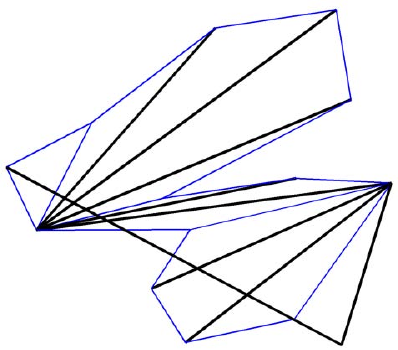}\\[10pt]
\includegraphics[width=0.7\textwidth]{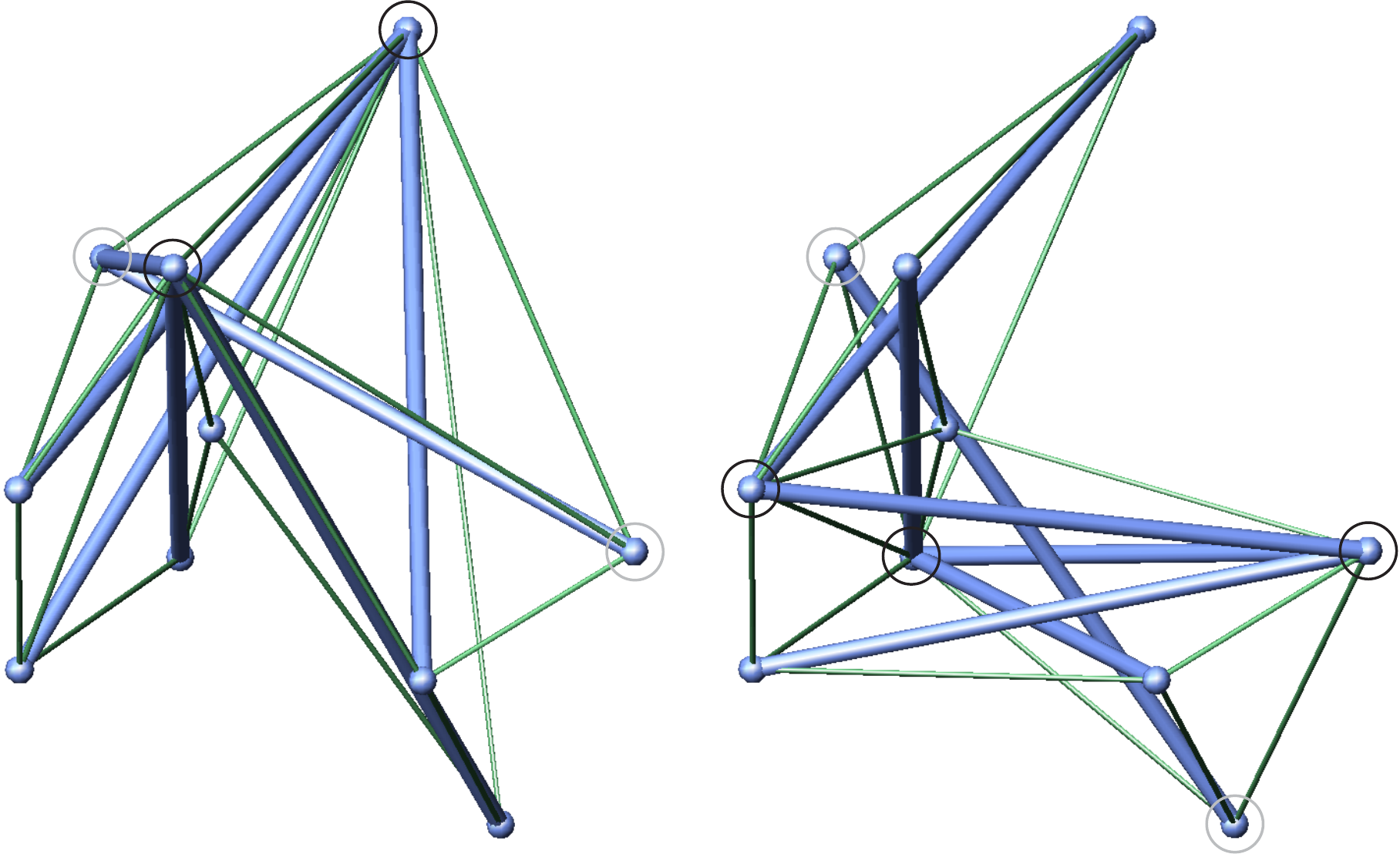}\,
\caption{Top: two universally rigid frameworks obtained from the same randomly generated set of nodes in $E^2$. Bottom: two constructions for the same randomly generated set of nodes in $E^3$. The frameworks on the right are composed by two fans. Thin and thick edges correspond respectively to positive and negative stresses. In the bottom pictures, black circles locate the end nodes of central edges, and grey circles locate the neighboring nodes. }
\label{ex1}
\end{center}
\end{figure}

\noindent Figure \ref{ex1} shows 2D and 3D examples obtained by applying the construction to a randomly generated set of nodes. The sign of the selfstress corresponding to a positive definite $\Omega$ is also represented.

Notice that we can find other generalized Gr\"unbaum frameworks. For example, the one shown in Fig.\,\ref{ex1} (top right) has two centers, corresponding to two fans with one side in common. It is easy to see that, in order for multi-fan frameworks to be universally rigid, the centers should be on opposite sides of the edge connecting the neighbors. Multi-fan frameworks exist also in three dimensions, like for example the one in Fig.\,\ref{ex1} (bottom right).

%

\section{Concluding remarks}

We have given a construction for generic universally rigid frameworks in two and three dimensions with the minimum number of edges, with a significant improvement over existing methods. Other than having applications to sensor networks, multi-agent systems and protein conformation analysis, the construction can also be used to design super-stable structural and mechanical systems with guaranteed strength and stiffness properties.

Our constructions relies on the computation of the convex hull of the set of nodes, a relatively quick operation, even if performed dynamically, i.e. if nodes are sequentially added and removed (see e.g. \cite{cgal:hs-ch3-00a}). Thus, in case of a sensor network or a multi-agent system with a variable number of sensors/agents, it would be easy to update the network of communication links in real, time while preserving the universal rigidity of the underlying framework.

As already mentioned, for nongeneric configurations the number of edges can be further reduced, with the limitation that any change of the configuration must happen on a lower-dimensional nongeneric manifold. The problem of constructing minimal nongeneric universally rigid frameworks will be the subject of future work.

\section*{Acknowledgements}

Support for this work was provided in part by the National Science Foundation grant CMMI-1000652.



\bibliographystyle{plain}
\bibliography{formation}

\end{document}